\documentclass[a4paper]{amsart}
\pdfoutput=1
\theoremstyle{remark}
\newtheorem*{example}{Example}
\newtheorem*{remarks}{Remarks}
\newtheorem*{acknowledgements}{Acknowledgements}
\newtheorem*{proofofstar}{Proof of $\eqref{eq:roots}$}

\newcommand{\field}[1]{\mathbb{#1}}
\newcommand{\Q}{\field{Q}}
\newcommand{\R}{\field{R}}
\newcommand{\Z}{\field{Z}}
\newcommand{\lra}{\longrightarrow}
\newcommand{\lms}{\longmapsto}
\newcommand{\wt}{\widetilde}
\newcommand{\sett}[2]{\{#1\,|\,#2\}}
\newcommand{\ovl}[1]{\overline{#1}}
\newcommand{\x}{{\bf x}}
\renewcommand{\t}{{\bf t}}
\newcommand{\s}{{\bf s}}
\DeclareMathOperator{\charak}{char}

\DeclareMathOperator{\Tr}{Tr}
\DeclareMathOperator{\Inn}{Inn}

\newsymbol\nmid      232D

\title{A twisted Laurent series ring that is a noncrossed product}
\date{\today}

\author{Timo Hanke{${}^1$}}
\address{
  University of California at San Diego\\
  Department of Mathematics, 0112\\
  9500 Gilman Dr.\\
  San Diego, CA 92093-0112\\
  USA \\
  e-mail : thanke@math.ucsd.edu
}

\begin{document}


\begin{abstract} 
The striking results on noncrossed products
were their existence (Amitsur)
and the determination of $\Q(t)$ and $\Q((t))$ as their smallest possible centres
(Brussel).
This paper gives the first fully explicit
noncrossed product example over $\Q((t))$.
As a consequence,
the use of deep number theoretic theorems
(local-global principles such as the Hasse norm theorem
and density theorems)  
in order to prove existence is eliminated.
Instead, the example can be verified by direct calculations.
The noncrossed product proof is short and elementary.
\end{abstract}

\maketitle

\begin{quotation}
  \tiny NOTICE: this is the corrected and updated author's version of a work that was
previously published in Israel Journal of Mathematics. Changes resulting from the publishing
process, such as peer review, editing, corrections, structural formatting, and
other quality control mechanisms may not be reflected in this document. 
On the other hand, changes
have been made to this work since the original publication 
in Israel Journal of Mathematics, vol.\ 150 (2005). 
\end{quotation}

\footnotetext[1]{Supported in part by the DAAD (Kennziffer D/02/00701).}

Let $K$ be the maximal real subfield of the $7$-th cyclotomic field,
the cubic Galois number field with discriminant $49$.
We define a cubic cyclic division algebra $D$ over $K$ 
and a $\Q$-automorphism $\sigma$ of $D$, 
such that the twisted Laurent series ring $D((\x,\sigma))$ 
is a noncrossed product of index and exponent $9$.  

Choose the primitive element $\alpha := \zeta+\zeta^{-1}\in K$
where $\zeta$ denotes a primitive \mbox{$7$-th} root of unity.
A non-trivial automorphism of $K$ over $\Q$ is defined by
$\sigma:\alpha\lms-\alpha^2-\alpha+1$.
The minimal polynomial of $\alpha$ is $x^3+x^2-2x-1$.
We immediately observe that the prime $2$ is unramified in $K$ 
because this polynomial is irreducible modulo $2$.
The prime $7$ is ramified since $\pi:=\alpha^2-\alpha-2$ is 
an element with $N_{K/\Q}(\pi)=7$.
Both observations about ramification are also clear from the general theory of cyclotomic fields.

Let $L=K(\theta)$ be the cubic extension of $K$ 
generated by a root $\theta$ of the polynomial
$$
f(x) = x^3+(\alpha-2)x^2-(\alpha+1)x+1\in K[x].
$$
Then $L/K$ is cyclic and a non-trivial automorphism of $L/K$ is defined by
$\phi :\theta\lms-\theta^2+(-\alpha+1)\theta+2$.
Let $D$ be the cyclic cubic algebra $(L/K,\phi,2\pi)$, i.e.
$$
D=L\oplus Lu \oplus Lu^2, \quad 
u^3=2\pi, \quad
ul=\phi(l)u \quad \forall l\in L.
$$ 
It is easily verified that $f$ has no roots modulo the prime ideals
$(2)$ and $(\pi)$ of $K$,
hence $(2)$ and $(\pi)$ are inertial in $L/K$.
Therefore, 
$D$ and its completions at the primes $(2)$ and $(\pi)$ of $K$ are division algebras.
The automorphism $\sigma$ of $K$ extends to an outer automorphism $\wt\sigma$ of $D$
by defining
$$
\wt\sigma(\theta) := \frac{1}{673}\sum_{i,j=0}^{2}c_{ij}\theta^iu^j,
\qquad \wt\sigma(u) := \lambda u,
$$
where
\begin{align*}
c_{ij} &:= 
\begin{bmatrix}
  303\alpha^2 - 154\alpha - 276 &  314\alpha^2 + 218\alpha - 326 &  -48\alpha^2 + 151\alpha + 157\\
  390\alpha^2 + 708\alpha - 855 &   40\alpha^2 - 238\alpha + 430 & -397\alpha^2 -  27\alpha + 275\\
 -106\alpha^2 +  25\alpha + 543 & -128\alpha^2 -  46\alpha -  30 &  135\alpha^2 +  38\alpha -  63\\
\end{bmatrix}, \\
\lambda &:= (\alpha^2+\alpha)+(-\alpha+1)\theta-\theta^2. 
\end{align*}
The relations to verify at this point are
$$
\wt\sigma(u)^3 = \sigma(2\pi), \quad
\wt\sigma(u)\wt\sigma(\theta) = \wt\sigma\phi(\theta)\wt\sigma(u), \quad 
f^\sigma(\wt\sigma(\theta)) =0,
$$
where $f^\sigma$ is the polynomial obtained from $f$ by applying $\sigma$ to the coefficients.

\begin{remarks}
The images $\wt\sigma(\theta)$ and $\wt\sigma(u)$ were obtained as the first application of
a project on computational methods for simple algebras over number fields,
whose results are to be published in detail in separate papers.
The computation for the present case involves two relative norm equations,
one for each of the given images,
in fields of absolute degree $27$ and $9$ respectively.
The equation for $\wt\sigma(u)$ is rather obvious,
it is $N_{L/K}(\lambda) = \frac{\sigma(\pi)}{\pi}$.
All computations were carried out using the computer algebra system MAGMA \cite{magma}.
\end{remarks}

The twisted Laurent series ring $D((\x,\wt\sigma))$ 
is the ring of formal series
$\sum_{i\geq k} a_i \x^i$, $k\in\Z$, $a_i\in D,$
with addition defined componentwise and multiplication defined by the rule 
$\x a=\wt\sigma(a)\x$ for all $a\in D$.
It is a division algebra of degree $9$.
The centre consists of the series $\sum_{i\geq k}c_i\t^i$, $c_i\in\Q$,
with $\t=d^{-1}\x^3$ 
for an element $d\in D^\times$ with
$\wt\sigma^3=\Inn(d)$ and $\wt\sigma(d)=d$.
Since $\t$ is a commutative indeterminate over $\Q$,
the centre can be identified with the Laurent series field $\Q((\t))$.
For a reference on the facts stated in this paragraph see Pierce \cite[\S~19.7]{pierce:ass-alg}.
An element $d$ as above is, for instance, 
$d := \sum_{i,j=0}^{2}d_{ij}\theta^iu^j$ with
$$
d_{ij} :=  
\begin{bmatrix}
 468\alpha^2+536\alpha+136  & 52\alpha^2+30\alpha-184   & -126\alpha^2-357\alpha+77 \\
 628\alpha^2-244\alpha-624  & -350\alpha^2+14\alpha+574 & 151\alpha^2-163\alpha-201 \\ 
 -416\alpha^2-240\alpha+324 & 14\alpha^2-14\alpha+84    & 74\alpha^2+166\alpha-124  \\
\end{bmatrix}.
$$

The remaining part of the paper shows that $D((\x,\wt\sigma))$ is a noncrossed product.
A valuation on a division ring $D$ is a map $v: D\lra\R$ with
\begin{align*}
&v(x)\geq 0, \quad\text{and}\quad v(x)=0 \iff x=0, \\
&v(xy)=v(x)v(y), \\
&v(x+y) \leq \max\{v(x),v(y)\}, 
\end{align*}
for all $x,y\in D$.
Associated to $v$ is
the {\em valuation ring}
$B_v := \sett{x\in D}{v(x)\leq 1}$
with unique maximal ideal
$M_v := \sett{x\in D}{v(x)< 1}$,
and the {\em residue division ring}
$\ovl D:=B_v/M_v$. 
If $L \subseteq D$ is a (commutative) subfield
then $v|_L$ clearly is a non-archimedian valuation on $L$.
Moreover, $\ovl L$ is a subfield of $\ovl D$.
If $Z(D)\subseteq L$ then $v|_L$ is the unique extension of $v|_{Z(D)}$ 
to a valuation on $L$.
References for valuations on division rings are Wadsworth \cite{wadsworth:survey}
and Schilling \cite{schilling:th-of-val}.

\begin{example}
  \label{ex:x-adic}
A valuation on $D((\x,\wt\sigma))$, called the $\x$-adic valuation, is defined by
$$ 
v(\sum_{i\geq k} a_i\x^i) := \delta^{\min\sett{i\in\Z}{a_i\neq 0}}, 
\quad 0<\delta<1.
$$
We have
$ B_v = \sett{\sum_{i\geq 0} a_i \x^i}{a_i\in D}$ and
$M_v = \sett{\sum_{i\geq 1} a_i \x^i}{a_i\in D}$.
The residue ring $\ovl{D((\x,\wt\sigma))}$ is routinely identified with $D\x^0$.
$v|_{\Q((\t))}$ is the $\t$-adic valuation.
\end{example}

Write $\mu_n\subset F$ if the field $F$ contains a primitive $n$-th root of unity.
We will make use of a standard argument in valuation theory~:
\begin{equation}
\label{eq:roots}
\text{\begin{minipage}{9.5cm}If a totally and tamely ramified finite extension 
of discretely valued fields  $E/F$ is cyclic then $\mu_n\subset\ovl F$
where $n=[E:F]$.
\end{minipage}}
\end{equation}
A short proof of this fact is included for completeness at the end of this paper.

Assume there is a maximal subfield $M$ of $D((\x,\wt\sigma))$
that is Galois over the centre $\Q((\t))$.
Then $[M:\Q((\t))]=9$. 
We claim that $M/\Q((\t))$ is unramified.
For, if it is totally ramified then \eqref{eq:roots} yields the contradiction $\mu_9\subset\Q$.
If the ramification index is $3$ then \eqref{eq:roots} applied to $M$ 
over the maximal unramified subextension 
states $\mu_3\subset\ovl M$.
This is a contradiction since $\ovl M$ is a cubic extension of $\Q$.
Therefore, $M/\Q((\t))$ is unramified.

It follows that $\ovl M/\Q$ is Galois of degree $9$.
Since the degree is $9$, $\ovl M$ is a maximal subfield of $D$ and contains $K$.
Above it was shown that the completions of $D$ at the primes $(2)$ and $(\pi)$ of $K$ are division algebras.
This means that $2$ and $7$ extend uniquely to primes of $\ovl M$.
If $\ovl M/\Q$ is cyclic then the subfields of $\ovl M$ are linearly ordered,
hence $K$ is the unique cubic subfield.
Since $7$ ramifies in $K$ it follows that 
the inertia field of $7$ in $\ovl M$ is $\Q$ itself.
This means that $7$ is totally ramified in $\ovl M/\Q$,
hence $\mu_9\in\ovl\Q_7$ by \eqref{eq:roots}, a contradiction.
If $\ovl M/\Q$ is not cyclic there exists a cubic cyclic subfield $K'$ of $\ovl M$ linearly disjoint to $K$.
Since the prime $2$ is inertial in $K$ it must be ramified in $K'$,
hence $\mu_3\subset\ovl\Q_2$ by \eqref{eq:roots}, a contradiction.
We have shown that $M$ as assumed cannot exist,
i.e.\ $D((\x,\wt\sigma))$ is indeed a noncrossed product.

\begin{remarks}
The exponent is $9$ because of the known fact that every division algebra over $\Q((\t))$ of prime exponent is a crossed product.

Analogously, the twisted function field $D(\x,\wt\sigma)$ is a division algebra of degree $9$ over $\Q(\t)$.
It is also a noncrossed product since $D((\x,\wt\sigma))$ is the completion of $D(\x,\wt\sigma)$ 
with respect to the $\x$-adic valuation.

The technique of the proof above is similar to a corresponding part in Brussel \cite{brussel:noncr-prod}.
The existence of noncrossed products of the form $D((\x,\wt\sigma))$ 
was shown before in \cite{hanke:thesis} using local-global principles.
However it was so far impossible to compute the outer automorphisms explicitly. 
Another explicit noncrossed product example over the centre $\Q((\s))((\t))$
(and $\Q(\s)(\t)$ resp.) can be found in \cite{hanke:expl-ex}.
\end{remarks}

\begin{proofofstar}
Let $E/F$ be a totally and tamely ramified finite cyclic extension of discretely valued fields of degree $n$.
Denote the Galois group by $G$ and the discrete valuation by $v$.
Since $v$ extends uniquely from $F$ to $E$, the elements
$\Delta_\sigma a:=\frac{\sigma(a)}{a}$ are valuation units for all 
$\sigma\in G$ and $a\in E^\times$.
Moreover, since $E/F$ is totally ramified, $\ovl{\tau\Delta_\sigma a}=\ovl{\Delta_\sigma a}$ for all $\tau\in G$.
It follows 
$\ovl{\Delta_{\sigma\tau}a}=
\ovl{\sigma\Delta_\tau a}\cdot\ovl{\Delta_\sigma a}=
\ovl{\Delta_\sigma a}\cdot\ovl{\Delta_\tau a}$
and, in particular, $\ovl{\Delta_\sigma a}^n=\ovl{\Delta_1 a}=\ovl{1}$.

Now let $\sigma$ be a generator of $G$ and let $a$ be a uniformizer of $E$. 
In order to show that $\ovl{\Delta_\sigma a}$ is a primitive $n$-th root of unity
we assume $\ovl{\Delta_\sigma a}^m=\ovl{1}$ and deduce $n|m$.
Setting $b:=a^m$, it is $\ovl{\Delta_\sigma b}=\ovl{1}$.
Since $\sigma$ generates $G$, $\ovl{\Delta_{\tau} b}=\ovl{1}$ follows for all $\tau\in G$.
By the hypothesis of tame ramification, $\charak\ovl{F}\nmid n$, hence
$ \sum_{\tau\in G} \ovl{\Delta_\tau b} =\ovl{n}\neq\ovl{0} $.
This means $v(\sum_{\tau\in G} \Delta_\tau b)=1$ so that
$$
v(a)^m=v(b)=v(b\sum_{\tau\in G}\Delta_\tau b)=
v(\sum_{\tau\in G}\tau(b))=v(\Tr_{E/F}(b))\in v(F).
$$
Since $E/F$ is discrete and totally ramified of degree $n$,
and $a$ was assumed to be a uniformizer of $E$,
it follows that $n|m$. \qed
\end{proofofstar}

\begin{acknowledgements}
I would like to thank the Department of Mathematics at UCSD and in particular my host Adrian Wadsworth
for their kind hospitality and support. 
The fellowship from the DAAD (German Academic Exchange Service) for this visit is greatly acknowledged.
Thanks also go to Joachim Gr\"ater for his support in previous work 
that contributed to this result.
Finally, I would like to express my gratitude to the MAGMA group for making available their computer algebra system 
and in particular to Claus Fieker for his helpful correspondence regarding the norm solver module.
\end{acknowledgements}

\bibliographystyle{amsplain}
\bibliography{thanke_2005_ijmath}

\nocite{amitsur:central-div-alg}
\nocite{brussel:noncr-prod}

\end{document}